\documentclass[12pt]{article}
\usepackage{color}
\usepackage{amsmath,amsfonts,amssymb,graphics,amsthm}
\usepackage {graphicx,microtype}
\usepackage[margin=10pt,font=small,labelfont=bf]{caption}
\usepackage[usenames,dvipsnames,svgnames,table]{xcolor}
\usepackage[colorlinks=true,pdfpagemode=none,citecolor=OliveGreen,linkcolor=BrickRed,urlcolor=BrickRed,pdfstartview=FitH]{hyperref}
\usepackage{sidecap}

\newtheorem{theorem}{Theorem}[section]
\newtheorem{corollary}[theorem]{Corollary}
\newtheorem{lemma}[theorem]{Lemma}

\newtheorem{conjecture}[theorem]{Conjecture}

\theoremstyle{remark}
\theoremstyle{remark}\newtheorem{remark}[theorem]{Remark}

\def\MR#1{\href{http://www.ams.org/mathscinet-getitem?mr=#1}{MR#1}}
\def\suml{\sum\limits}
\newcommand{\tree}{\mathcal T_o}
\newcommand{\trunk}{\mathrm{Ray}(o,r)}

\newcommand{\LE}{\mathbf{LE}}

\newcommand{\wt}{\widetilde}

\renewcommand{\path}{ \mathcal{P}}

\newcommand{\N}{\mathbb{N}}

\newcommand{\Z}{\mathbb{Z}}

\newcommand{\simC}{\sim_C}
\newcommand{\vtau}{T}
  
\DeclareMathOperator{\Pmathop}{\mathbb{P}}
\DeclareMathOperator{\Emathop}{\mathbb{E}}
\renewcommand\P{\Pmathop\mathopen{}}
\newcommand\E{\Emathop\mathopen{}}

\def\Psub_#1{\P_{\! #1}}

\def\Pbig#1{\P\mkern-.5mu\bigl[#1\bigr]}
\def\Psubbig_#1#2{\Psub_{#1}\mkern-1.5mu\bigl[#2\bigr]}

\def\Psubbigg_#1#2{\Psub_{#1}\mkern-1.5mu\biggl[#2\biggr]}

\def\Ebig#1{\E\mkern-1.5mu\bigl[#1\bigr]}
\def\Esubbig_#1#2{\E_{#1}\mkern-1.5mu\bigl[#2\bigr]}
\def\EsubBig_#1#2{\E_{#1}\mkern-1.5mu\Bigl[#2\Bigr]}
\def\Ebigg#1{\E\mkern-1.5mu\biggl[#1\biggr]}
\def\Esubbigg_#1#2{\E_{#1}\mkern-1.5mu\biggl[#2\biggr]}
\def\EBig#1{\E\mkern-1.5mu\Bigl[#1\Bigr]}
\def\vertspace{\mkern1.5mu}

\DeclareMathOperator{\WSF}{WSF}

\begin {document}
\title{Occupation measure of random walks\\
	and wired spanning forests\\ in balls of Cayley graphs}
\author{
	Russell Lyons\thanks{Department of Mathematics, Indiana
		University. Partially supported by the National
		Science Foundation under grant DMS-1007244 and by Microsoft Research.
		Email: \protect\url{rdlyons@indiana.edu}.}
	\and
	Yuval Peres\thanks{Microsoft Research, Redmond, WA.
		Email: \protect\url{peres@microsoft.com}.}
	\and
	Xin Sun\thanks{Department of Mathematics, Columbia University. Partially supported by
		Microsoft Research and by Simons Society of Fellows.
		Email: \protect\url{xinsun@math.columbia.edu}.}
	\and
	Tianyi Zheng\thanks{Department of Mathematics, UCSD.
	           Email: \protect\url{tzheng2@math.ucsd.edu}.}	
	           }
\date{}
\maketitle
\begin{abstract}
	We show that 
	for finite-range, symmetric random walks on general transient Cayley
   graphs, the expected occupation time of any given ball of radius $r$ is
   $O(r^{5/2})$. 
	We also study the volume-growth  property of the wired spanning forests on general Cayley graphs, showing that the expected number of vertices in the
	component of the identity inside any given ball of radius $r$ is
   $O(r^{11/2})$.
\end{abstract}

\begin{abstract}
	On montre que toute marche al\'eatoire sym\'etrique \`a pas born\'es sur un graphe de Cayley
   transitoire satisfait que l'esp\'erance du temps d'occupation d'une boule quelconque de rayon $r$
   vaut $O(r^{5/2})$.
   On \'etudie aussi la croissance du volume des f\^orets recouvrantes c\^abl\'ees dans les graphes
   de Cayley g\'en\'eraux, en montrant que l'esp\'erance du nombre de sommets appartenant \`a la
   composante connexe de l'identit\'e dans une boule quelconque de rayon
   $r$ vaut $O(r^{11/2})$.
\end{abstract}

\section{Introduction}
Given a transient, symmetric random walk $S$ starting from a vertex $o$ in a
Cayley graph $G=(V,E)$, let $ L_r:=\bigl|\bigl\{t:S_t\in
B(o,r)\bigr\}\bigr|$, where $B(o, r)$ is the set of vertices within graph
distance $r$ of $o$.
Suppose for the moment that $S$ is simple random walk.
If $G$ has polynomial growth of degree $D$, then Varopoulos' estimate $p_t(o, x) \lesssim t^{-D/2}$ (see, e.g., \cite[Corollary 7.3]{CGP})
yields $\E[L_r] \lesssim r^2$ (see Remark~\ref{rmk:polynomial}). 
Here, $a(t)\lesssim b(t)$ means that $\exists\, c>0$ such that $a(t)\le c\,b(t)$ for all $t$.
Similarly, Varopoulos' estimate $p_t(o, o) \lesssim e^{-c t^{1/3}}$ for groups of exponential growth (see \cite[Corollary 7.4]{CGP}) 
yields $\E[L_r] \lesssim r^3$ (see the proof of \cite[Proposition 2.3]{Cubic}). 
When the walk escapes at a linear rate, a simple argument (Lemma~\ref{lem:linear}) shows that $\E[L_r]\lesssim r$. In particular, the linear bound holds for nonamenable Cayley graphs.
We believe that the following  quadratic bound holds in general; to 
the best of our knowledge, this is open.
\begin{conjecture}\label{conj:quadratic}
	For a symmetric random walk $S$ on a transient Cayley graph $G$, let $L_r$ be the occupation time of $B(o,r)$ defined as above. Then $\E[L_r]\lesssim r^2$.
\end{conjecture}
As an example of amenable Cayley graphs of exponential growth where a
quadratic bound is easy to establish, consider simple random walk on
lamplighter groups over any base
group which has polynomial growth or, more generally, any base group known
to have quadratic occupation time: we can bound the occupation time of
balls in the Cayley graph of the lamplighter
group by the occupation time of balls of the projection of simple
random walk under the quotient map to the base group.
In this paper, although we cannot prove Conjecture~\ref{conj:quadratic}, we
establish a general $5/2$-power bound for finite-range, symmetric random walks
(i.e., symmetric random walks whose jumps have bounded support).
\begin{theorem}\label{thm:occupation-time}
	Let $G$ be a transient Cayley graph and $V(r):=|B(o,r)|$. Then for every finite-range,
	symmetric random walk on $G$,
	\begin{equation}\label{eq:occupation-thm}
	\E[L_r] \lesssim r^2\sqrt{\log V(r)}\,.
	\end{equation}In particular, $\E[L_r]\lesssim r^{5/2}$.
\end{theorem}

By comparison, if $\tau_r$ denotes the first exit time of $B(o,r)$ of a
symmetric random walk starting at $o$, it is known that $$E[\tau_r] \lesssim
r^2$$ for all Cayley graphs. (See Theorem~\ref{thm:escape} for a proof.)

Let $G=(V,E)$ be an infinite graph. The wired spanning forest measure on $G$ is defined as the infinite-volume limit of the wired spanning tree measures on a
sequence of finite subgraphs exhausting $G$:  Let $V_1\subset V_2\subset \dotsm $ be finite subsets of $V$ whose induced subgraphs $G_n$ are connected with
$\bigcup_{n=1}^\infty V_n=V$. Let $\mu^F_n$ be the uniform spanning tree measure on $G_n$. Then as
a probability measure on edge configurations, $\mu^F_n$ restricted to any finite subset of $E$ converges. This defines a unique probability measure
$\mu^F$ on $2^E$, which we call the \emph{free spanning forest}.  Another way of taking limits of spanning trees is as follows. Suppose $G_n$ are defined as
above. Let $G_n^W$ be obtained from $G$ by identifying all the vertices outside $G_n$ to one new vertex
and $\mu^W_n$ be the uniform spanning tree measure on $G^W_n$. Then
$\mu^W_n$ also has a limit $\mu^W$, which we call the \emph{wired spanning
forest}.  These results are due to \cite{Pemantle}.
The free and wired spanning forests are the same if $G$ is of polynomial
growth or, more generally, amenable \cite{Pemantle,benjamini2001special}. They can be different, such as on the Cayley graph of a free group.  See
\cite{benjamini2001special,TreeBook} for more details.

On Cayley graphs, the wired spanning forest (WSF) has a single component if the graph has at most quartic growth; otherwise, there are infinitely many components in
the WSF \cite{Pemantle}. In the latter case, the geometry of the WSF has intriguing behaviors. 
Let $\tree$ be the component containing $o$ in the WSF of $G$. 
For Cayley graphs with polynomial growth of order at least
quartic, $\Ebig{\vertspace|\tree\cap B(o,r)|\vertspace}\asymp r^4$, whereas
nonamenable Cayley graphs satisfy $\Ebig{\vertspace|\tree\cap B(o,r)|\vertspace} \asymp r^2$ \cite[Section
13]{benjamini2001special}. Here, $a(t)\asymp b(t)$ means $a(t)\lesssim b(t)$ and $b(t)\lesssim a(t)$.
In \cite{GeometryofUSF}, the authors provided a detailed analysis of
the geometry of the WSF on $\Z^d$ ($d\ge 5$). Among other results, they showed that the tree components have ``stochastic dimension" 4.  In this paper, we
extend the investigation of the volume-growth property of the WSF to general Cayley graphs (Theorems \ref{thm:whole-volum} and \ref{thm:C(0,r)}). 

Using a similar method as we use to prove Theorem \ref{thm:occupation-time}, we show the following upper bound:

\begin{theorem}\label{thm:whole-volum}
	Let $G$ be a Cayley graph and $V(r)=|B(o,r)|$. Then
	\begin{equation}\label{eq:whole-volum}
	\Ebig{\vertspace|\tree\cap B(o,r)|\vertspace}\lesssim
   r^{4}\log^{3/2} V(r)\,.
	\end{equation}In particular, $\Ebig{\vertspace|\tree\cap
   B(o,r)|\vertspace}\lesssim r^{11/2}$.
\end{theorem}

Let $C(o,r)$ be  the connected component of $\tree\cap B(o,r)$ containing $o$.
This provides another way to measure the growth of the WSF.
We show the following upper bound in terms of the exit time $\tau_r$ for random walk:

\begin{theorem}\label{thm:C(0,r)}
	Given a Cayley graph $G$  of superpolynomial growth,  let $C(o,r )$ be defined as above. Then there exists $r_0$ such that
	\begin{equation}\label{eq:C(o,r)}
	\Ebig{\vertspace|C(o,r)|\vertspace} \le 4\E^2[\tau_{6r}]\quad \textrm{for}\quad r>r_0\,.
	\end{equation} 
\end{theorem}
\begin{remark}
	As will be clear from our proof of \eqref{eq:C(o,r)}, the constants involved are not optimal.
\end{remark}

For Cayley graphs of polynomial growth, we have $\Ebig{\vertspace|C(o,r)|\vertspace}\le \Ebig{\vertspace|\tree\cap B(o,r)|\vertspace}\lesssim r^{4}$. Since $\E[\tau_r] \lesssim r^2$
for all Cayley graphs, Theorem~\ref{thm:C(0,r)} implies that $\Ebig{\vertspace|C(o,r)|\vertspace}
\lesssim r^4$ in general. We believe that $\Ebig{\vertspace|\tree\cap B(o,r)|\vertspace} \asymp
\Ebig{\vertspace|C(o,r)|\vertspace}$  and
hence $\Ebig{\vertspace|\tree\cap B(o,r)|\vertspace} \lesssim r^4$ for general Cayley graphs.

\bigskip
\noindent{\bf Acknowledgments.}
We are grateful to Terry Tao for providing the reference \cite{GreenTaoB}.
We thank the referees for useful comments.
This work was begun while the third author was an intern in the Theory Group
at Microsoft Research, Redmond.

\section{Occupation measure of random walks}
\subsection{Preliminaries}\label{subsec:pre}

The only random walks $S = (S_0, S_1, \ldots)$ on groups that we consider
are those where for all $t \ge 1$, the random variables $S_{t-1}^{-1} S_t$
are independent and identically distributed.
Such a random walk is called \emph{symmetric} if for all $g$, we have
$\Pbig{S_0^{-1} S_1 = g} = \Pbig{S_0^{-1} S_1 = g^{-1}}$.
We usually choose $S_0$ to be the identity, $o$.

Suppose $\Gamma$ is a group generated by a finite subset $X$, i.e., every element in $\Gamma$ can be written as a product of elements  in $X\cup X^{-1}$.
The Cayley graph $G$ associated to $(\Gamma,X)$ is the unoriented graph with vertices
$\Gamma$ and edges $\bigl\{[g,gx]:g\in \Gamma,\,x\in X\bigr\}$. Every Cayley graph is a
connected, vertex-transitive graph. 

For a Cayley graph $G$, a vertex $o\in G$, and $r>0$, let $d_G$ denote the graph distance in $G$ and  $B(o,r):=\{v\in G:d_G(o,v)\le r\}$. We call
$V(r):=|B(o,r)|$ the \emph{volume function} of $G$. Due to Gromov's theorem \cite{Gromov:poly},
it is well known that either $V(r)\asymp r^D$ for some $D\in \N$ or $\lim_{r\to\infty} V(r)/r^D=\infty$ for
all $D\in \N$. In the
former case, we say that $G$ has \emph{polynomial growth of degree $D$}. In
the latter case, we say that $G$ has \emph{superpolynomial growth}.
These properties are independent of the choice of the generating set $X$ of $G$. 

Given a Cayley graph $G$ with $d := |X\cup X^{-1}|$,
lazy simple random walk on $G$ is the Markov chain $S=(S_t)^\infty_{t=0}$ on $\Gamma$ with transition probabilities $p(g, gx)
= 1/(2d)$ for $x \in X \cup X^{-1}$ and $p(g, g) = 1/2$.
We assume that the identity is not an element of $X$.

The following facts concerning the occupation time $L_r$ and the escape time $\tau_r$  are not needed for the rest of the paper. We record them for completeness. 

\begin{lemma}\label{lem:linear}
	Suppose $S_t$ is a random walk on a Cayley graph $G$ such that
	$\liminf_{t \to\infty} d_G(o, S_t)/t > 0$ a.s. Then  $\E[L_r]\lesssim r$.
\end{lemma}
\begin{proof}
	We may choose $\epsilon > 0$ and $t_0 < \infty$ 
	so that 
	\[
	\Psub_o[\forall t \ge t_0 \quad d_G(o, S_t) > \epsilon t] >
	1/2\,.
	\]
	Let $s(r):= \max \{ 2r/\epsilon, t_0\}$.
	Then for every $t$, we have 
	\[
	\Pbig{\forall m \ge s(r)\quad
		S_{t + m} \notin B(o, r) \bigm| S_t \in B(o, r)} > 1/2\,,
	\]
	so $\E[L_r] < 2s(r)$.
\end{proof}

Note that if $\Gamma$ is a nonamenable group, then the hypothesis of Lemma
\ref{lem:linear} holds: \cite{Kesten1,Kesten2} showed that there is some
$\rho < 1$ such that for all $x \in \Gamma$ and all $t \in \N$,
we have $p_t(o, x) \le \rho^t$. The result then follows from a
Borel--Cantelli argument.

The following argument was
noted by Anna Erschler (personal communication, 2005). 

\begin{theorem}\label{thm:escape}
	$\E[\tau_r] \lesssim r^2$ for symmetric random walks on Cayley graphs.
\end{theorem}
\begin{proof}
	
	Because of the linear bound on nonamenable Cayley graphs even for occupation
	time (Lemma~\ref{lem:linear}) and of the stochastic domination of $\tau_r$ by $L_r$,
	it remains to show this bound on escape time when $G$ is amenable.
	Furthermore, we may assume that the support of the random walk generates
	the group $\Gamma$, as otherwise we take the subgroup it generates. 
	Let $W$ be a finite subset of the support of $S_1$ such that $W$ generates
	$\Gamma$. Because distances in any Cayley graph of $G$ differ from those in the Cayley graph
	generated by $W$ by a bounded factor, we may assume that $G$ is in fact
	the Cayley graph determined by $W$.
	We may also assume that the support of $S_1$ is contained in $B(o, 2r)$
	since if not, we may replace all jumps outside that ball by staying in
	place; the new random walk leaves $B(o, r)$ no earlier than the original
	random walk does.
	By \cite{Mok,KorSch}, there is a harmonic, equivariant, Hilbert-space valued,
	nonconstant function $\phi$ on $V$ (also see \cite[Theorem 3.1]{Lee-Peres}
	for an explicit construction).  Here, ``equivariant" means with respect to some
	affine isometric
	action of the group on the Hilbert space.
	Let $c := \Ebig{\|\phi(S_1) - \phi(o)\|^2} > 0$.
	\vadjust{\kern2pt}%
	Let $p_* := \min \{ p(o, x) : x \in W\}$.
	Then 
	\vadjust{\kern2pt}%
	$\|\phi(x) - \phi(y)\|^2 \le c/p_*$ when $x$ and $y$ are neighbors in $G$,
	\vadjust{\kern2pt}%
	whence $\|\phi(x) - \phi(y)\| \le \sqrt{c/p_*} \cdot d_G(x, y)$ for all vertices
	$x, y$ of $G$.
	\vadjust{\kern2pt}%
	In particular,
	$\|\phi(x) - \phi(o)\| \le 3r \sqrt{c/p_*}$ for $x \in B(o, 3r)$. 
	Since $\phi$ is harmonic, the sequence of random variables $\|\phi(S_n) - \phi(o)\|^2 - c n$
	\vadjust{\kern2pt}%
	forms a martingale, thus the optional-stopping theorem gives $\E\|\phi(S_{\tau_r}) - \phi(o)\|^2 = c\E[\tau_r]$. Since the support of $S_1$ is within $B(o,2r)$ and $\tau_r$ is the exit time of $B(o,r)$, the triangle inequality gives $S_{\tau_r}\in B(o,3r)$.
	Therefore
	\[\E[\tau_r] \le \big(3r\sqrt{c/p^*}\,\big)^2\cdot c^{-1}=9 r^2/p_*\,.\qedhere\]
\end{proof}

When the random walk
has bounded jumps, a stronger result
on the distribution of $\tau_r$ 
follows from the main result of
\cite{Lee-Peres-Smart}.

\subsection{Proof of Theorem~\ref{thm:occupation-time}}\label{subsec:proof}
There are three main ingredients in our proof of  Theorem \ref{thm:occupation-time}.
The first ingredient is a bound for the return probability of lazy random walks using the volume function $V(r)$, which is obtained in \cite{lyons2012sharp} by spectral embedding:
\begin{lemma}\label{lem:returning-estimate}
	Given a vertex-transitive graph $G$, let $p_m(o,o):=\P[S_m=o]$ be the
	return probability of a lazy, finite-range, symmetric random walk, $S$. Let $V$ be the volume function defined as above. Then there exist
	constants $c \in (0, 1)$ and $c' < \infty$ such that
	\begin{equation}\label{eq:returning-estimate}
	\forall m\in \N^+	\qquad p_m(o,o)\le c'm\int_0^1 \frac{e^{-\lambda m}}{V\bigl(c/\sqrt{\lambda}\,\bigr)} \,d\lambda \,.
	\end{equation}
\end{lemma}
\begin{proof}
	Combine Lemma 3.5 and Theorem 6.1 in \cite{lyons2012sharp}.
\end{proof}

The second ingredient is immediate from the main result of \cite
{Lee-Peres} in the amenable case and Lemma \ref{lem:linear} in the
nonamenable case:

\begin{lemma}\label{lem:diffusive}
	Given a vertex-transitive graph $G$, let $p_m(o,o):=\P[S_m=o]$ be the
	return probability of a lazy, finite-range, symmetric random walk, $S$.
   Then there exists a
	constant $c < \infty$ such that
	\begin{equation*}
	\forall r,n\in \N^+	\qquad \sum_{m=0}^n \Pbig{S_m \in B(o, r)} \le c r
   \sqrt n\,.
	\end{equation*}
\end{lemma}

The third ingredient is an important growth property of the volume function of Cayley graphs of superpolynomial growth, established in \cite{GreenTaoB}:
\begin{lemma}\label{lem:volume-function}
	Suppose $G$ is a Cayley graph of  superpolynomial growth. Then for all $k\in \N$, there exists $c_k>0$ such that
	\begin{equation}\label{eq:volume-function}
	\textrm{for all}\; a\ge 1\;\textrm{and}\;r\ge 1,
	\qquad \frac{V(ar)}{V(r)} \ge c_k a^k\,.
	\end{equation}
\end{lemma}
\begin{proof}
	This is an immediate consequence of \cite[Corollary 11.2]{GreenTaoB}.
\end{proof}

\begin{corollary}\label{cor:polybd}
	Suppose $G$ is a Cayley graph of superpolynomial growth of a group,
	$\Gamma$. Let $S$ be a lazy, finite-range, symmetric random walk on $G$ whose support
	generates $\Gamma$.  Write $p_m(x,y):=\Psub_x[S_m=y]$.
	Then there is a constant $c > 0$ such that for all $k\in \N^+$,
	there is some $c'' > 0$ (depending on $k$) such that
	for all $r, m\in \N^+$ and all $x, y \in \Gamma$,
	\begin{equation}\label{eq:polybd}
	p_m(x, y) \le c'' \bigl(m^{-k/2} r^k / V(r) + e^{-c^2 m/r^2}\bigr)
	\,.
	\end{equation}
\end{corollary}

\begin{proof}
	Choose $c$ as in \eqref{eq:returning-estimate}. 
	From the preceding two lemmas, we have
	\begin{align*}
	p_m(x, y)
	\le
	p_m(o, o)
	&\lesssim
	m\int_0^{1} \frac{e^{-\lambda m}}{V\bigl(c/\sqrt{\lambda}\,\bigr)}\, d\lambda 
	\\ &=
	m\int_0^{c^2/r^2} \frac{e^{-\lambda m}}{V\bigl(c/\sqrt{\lambda}\,\bigr)}\,
	d\lambda
	+ m\int_{c^2/r^2}^{1} \frac{e^{-\lambda m}}{V\bigl(c/\sqrt{\lambda}\,\bigr)}\, d\lambda\\
	&\lesssim
	\frac{m}{V(r)}\int_0^{c^2/r^2} \lambda^{k/2}r^{k} e^{-\lambda m}\,  d\lambda +
	m\int_{c^2/r^2}^{1}e^{-\lambda m}\, d\lambda\\
	&\lesssim
	m^{-k/2}r^{k}/V(r) +e^{-{c^2m}/{r^2}}\,,
	\end{align*}
	where in the last line, we use the change of variable
	$u:=m\lambda$. The implied constants depend on $k$.
	This proves \eqref{eq:polybd}.
\end{proof}

\begin{proof}[Proof of Theorem~\ref{thm:occupation-time}]
	We may clearly assume that the
	support of the walk generates the group, as otherwise we simply take the
	subgroup it generates together with a Cayley graph of the subgroup.
	We may also assume that $S$ is lazy, i.e., $p_1(o, o) \ge 1/2$.
	We wish to show that
	\begin{equation}\label{eq:occupation-time}
	\E[L_r] = \suml_{m=0}^\infty \Pbig{S_m\in B(o,r)}\lesssim
   r^2\sqrt{\log V(r)}\,.
	\end{equation}	
	Since the result is known for groups of polynomial growth, we assume $G$ is
	of superpolynomial growth.
	Write $\varphi(m)$ for the right-hand side of \eqref{eq:polybd}.
	Then
	$\forall m\in \N$ and $ r>0$,
	\begin{equation*}
	\Pbig{S_{m}\in B(o,r)} \le \varphi(m)V(r)\,.
	\end{equation*}
	Set $\alpha:=c^{-2}$, where $c$ is as defined in \eqref{eq:polybd}. Put
	\begin{align*}\label{eq:bound-returnning2}
	\Sigma^{(1)}_r&:=\suml_{m=0}^{\lfloor\alpha r^{2}\log V(r)\rfloor}\Pbig{S_m\in B(o,r)}\\
	\noalign{and}
	\Sigma^{(2)}_r&:=\suml_{ m> \alpha r^{2}\log V(r)}\varphi(m)V(r) 
	\,.
	\end{align*}
	
	By Lemma~\ref{lem:diffusive},
   \[
   \Sigma^{(1)}_r\lesssim r^{2}\sqrt{\log V(r)}\,.
   \]
   
   Since
   \[
   \suml_{m=0}^\infty \Pbig{S_m\in
		B(o,r)}\lesssim \Sigma^{(1)}_r+\Sigma^{(2)}_r\,,
   \]
   to prove \eqref{eq:occupation-time}, it suffices to show  that $\Sigma^{(2)}_r\lesssim r^2$.
	Choose $k > 2$ with Corollary \ref{cor:polybd} in mind.
	Now
	\begin{align}
	\suml_{ m> \alpha r^{2}\log V(r)} m^{-k/2}r^{k} \lesssim \bigl(r^{2}\log V(r)\bigr)^{-k/2+1}r^{k}\lesssim r^2 \,.
	\end{align}
	On the other hand,
	\begin{align}
	\suml_{ m> \alpha r^{2}\log V(r)} V(r)e^{-{c^2m}/{r^2}} \lesssim V(r)r^2e^{- \alpha c^2\log V(r)}=r^2\,.
	\end{align}
	Therefore, $\Sigma^{(2)}_r\lesssim r^2 $, as claimed.
\end{proof}

\begin{remark}\label{rmk:polynomial}
	If $G$ has polynomial growth, then we can separate the sum in
	\eqref{eq:occupation-time} at $\alpha r^2$ instead of at $\alpha r^2\log
	V(r)$. The same argument as above combined with the bounds $V(r) \asymp r^D
	$ and $p_{2m}(o,o)\asymp m^{-D/2}$ then gives a proof of the quadratic
   bound on occupation time; one does not need
   Lemma~\ref{lem:diffusive}, but only the trivial bound that every
   probability is at most 1.
\end{remark}

\section{Volume growth of the WSF}
Given a finite path $\path=\langle v_0,v_1, \dots,v_n\rangle $ in a graph $G$, we define the forward loop erasure of $\path$ (denoted by $\LE[\path]$)
by erasing cycles in $\path $ chronologically. More precisely, $\LE[\path]$ is defined inductively as follows. The first vertex $u_0$ of $\LE[\path]$ is
the vertex $v_0$ of $\path$. Supposing that $u_j$ has been set, let $k$ be the last index such that $v_k=u_j$. Set $u_{j+1}:=v_{k+1}$ if $k<n$; otherwise, let
$\LE[\path] :=\langle u_0,\dots, u_j \rangle$.  If $S$ is a simple random walk on a Cayley graph $G$,  then $\LE[S]$ is
called the \emph{loop-erased random walk (LERW)}.  There is no trouble defining the
forward loop erasure of $S$ a.s.\ if $G$ is transient.  For recurrent Cayley
graphs of quadratic growth, loop-erased random walk can be defined by
taking a limit (see
\cite{LawlerIntersections,benjamini2001special}). We omit the details, because we focus exclusively on transient graphs in the rest of the paper.

In \cite{WilsonAlgorithm}, Wilson discovered an algorithm for sampling uniform spanning trees on finite graphs using loop-erased random walk. In
\cite{benjamini2001special}, Wilson's algorithm was adapted to sample the WSF on transient graphs: Order the vertex set $V$ as $V=(v_1,v_2,\dots )$. Set $\mathcal{T}_0:=\emptyset$. Inductively, for each $n=1,2,\dots$,
run an independent simple random walk starting at $v_n$. Stop the walk when it hits $\mathcal{T}_{n-1}$ if it does; otherwise, let it run
indefinitely. Denote the resulting path by $\path_n$, and set $\mathcal T_n:=\mathcal T_{n-1}\cup \LE[\path_n]$. According to \cite[Theorem
5.1]{benjamini2001special}  no matter the ordering of $V$, the resulting forest is always distributed as the WSF on $G$. This method
of generating the WSF is called \emph{Wilson's method rooted at infinity}.

In fact, the theory of wired spanning forests extends to general networks,
i.e., general reversible random walks; see  \cite{benjamini2001special} or
\cite{TreeBook} for details. Thus, we will prove the following extension of
Theorem \ref{thm:whole-volum}:

\begin{theorem}\label{thm:gen-whole-volum}
	Let $G$ be a Cayley graph of a group $\Gamma$
	and $V(r):=|B(o,r)|$. Consider the WSF on $\Gamma$ corresponding to a finite-range symmetric random walk $S$ whose support generates $\Gamma$. Then
	\begin{equation}\label{eq:gen-whole-volum}
	\Ebig{\vertspace|\tree\cap B(o,r)|\vertspace}\lesssim
   r^{4}\log^{3/2} V(r)\,.
	\end{equation}In particular, $\Ebig{\vertspace|\tree\cap
   B(o,r)|\vertspace}\lesssim r^{11/2}$.
\end{theorem}

\begin{proof}
	The polynomial-growth case is known when the WSF is generated by simple random
	walk; the proof of its extension to finite-range symmetric random walks
	will be clear following Remark \ref{rmk:polynomial}.
	Thus, we assume $G$ has superpolynomial growth.  We may further assume that $S$ is lazy, since adding laziness simply produces loops in the random walk
	paths, which are then erased.
	
	Let $\{S^v\}_{v\in G}$ be a family of independent random walks with the
	same increment distribution as $S$ but
	such that $S^v$ starts from $v$. Let  $\Psub_v$ be the law of
	$S^v$. By Wilson's algorithm rooted at infinity,
	\begin{align}
	\nonumber
	\P[x\in \tree ]&\le \P[\exists y\in G\ \ \exists m\ge k\ge 0 \quad S^o(k)= S^x(m-k)=y]\\
	\label{eq:Wilson}
	&\le\suml_{y\in G}\suml_{m=0}^\infty\suml_{k=0}^m \Psub_o[S_k=y]\Psub_x[S_{m-k}=y]\,.
	\end{align}
	By reversibility and the Markov property,
	\begin{align*}
	\suml_{y\in G}\Psub_o[S_k=y]\Psub_x[S_{m-k}=y] = \Psub_o[S_m=x]\,.
	\end{align*}
	Combined with \eqref{eq:Wilson}, this leads to
	\begin{align*}
	\P[x\in \tree ]\le\suml_{m=0}^\infty(m+1)\Psub_o[S_m=x]\,.
	\end{align*}	
	Summing over $x\in B(o,r)$, we arrive at
	\begin{equation*}
	\Ebig{\vertspace|\tree \cap
		B(o,r)|\vertspace}\le\suml_{m=0}^\infty(m+1)\Psubbig_o{S_m\in B(o,r)}\,.
	\end{equation*}
	
	Decomposing this last sum similarly to the proof of Theorem \ref{thm:occupation-time}, we have
	\[
	\suml_{m=0}^\infty(m+1)\P[S_m\in B(o,r)]\lesssim \Sigma^{(3)}_r+\Sigma^{(4)}_r
	\,,
	\]where
	\begin{align*}
	\Sigma^{(3)}_r&:=\suml_{m=0}^{\lfloor\alpha r^2\log V(r)\rfloor}
	(m+1)\Psubbig_o{S_m\in B(o,r)}\,,\\
	\Sigma^{(4)}_r&:=\suml_{ m>\alpha r^2\log V(r)} V(r)(m+1) \varphi(m) 
	\,,
	\end{align*}
	and $\varphi$ is the right-hand side of \eqref{eq:polybd}.
	Using a very similar argument as in Theorem \ref{thm:occupation-time}, by
	choosing $k>4$ and  $\alpha:=2c^{-2}$, we obtain
   $$\Sigma^{(3)}_r\lesssim r^4\log^{3/2} V(r) \quad\textrm{and}\quad\Sigma^{(4)}_r\lesssim r^4\,,$$ thus concluding the proof. \qedhere
\end{proof}

To prove Theorem \ref{thm:C(0,r)}, we first  record an elementary fact concerning simple random walk on Cayley graphs.
\begin{lemma}\label{lemma:hitting-probablity}
	Let $G$ be a Cayley graph of superpolynomial growth and $S$ be a simple random walk starting from $o\in G$. For a vertex $x\in G$,
	let $|x|$ denote the graph distance from $x$ to $o$. Then for every $D>0$ there exists a positive constant $c_D$ such that
	\begin{equation}\label{eq:hitting-probability}
	\Psub_o[S\; \mathrm{hits}\; x]\le \frac{c_D}{|x|^D}\,.
	\end{equation}
\end{lemma}

\begin{proof}
Indeed, by Lemma~\ref{lem:returning-estimate}, for example,
	\[
	\Psub_o[S\; \mathrm{hits}\; x]\le
	\sum_{n \ge |x|} p_n(o, x)
	\lesssim \sum_{n \ge |x|} n^{-D-1} \lesssim |x|^{-D}\,.
	\qedhere
	\]
\end{proof}

\begin{proof}[Proof of Theorem \ref{thm:C(0,r)}]
	Suppose the WSF is generated via Wilson's algorithm by first sampling a
	simple random walk $S$ from $o$ and then sampling simple random walks from
	other vertices in a certain order. Let  $\mathrm{Ray}_o :=\LE[S]$ be the
	infinite ray emanating from $o$ in the WSF, $\>\mathrm{Ray}(o,r)
	:=\mathrm{Ray}_o\cap C(o,r)$, and $N_r :=|\mathrm{Ray}(o,r)|$. We first
	claim that $\E[N_r]\le  2\E[\tau_{3r}]$ for $r$ large enough.
	
	\begin{SCfigure}[50][ht!] \hspace{-.1\textwidth}
		\includegraphics[scale=0.5]{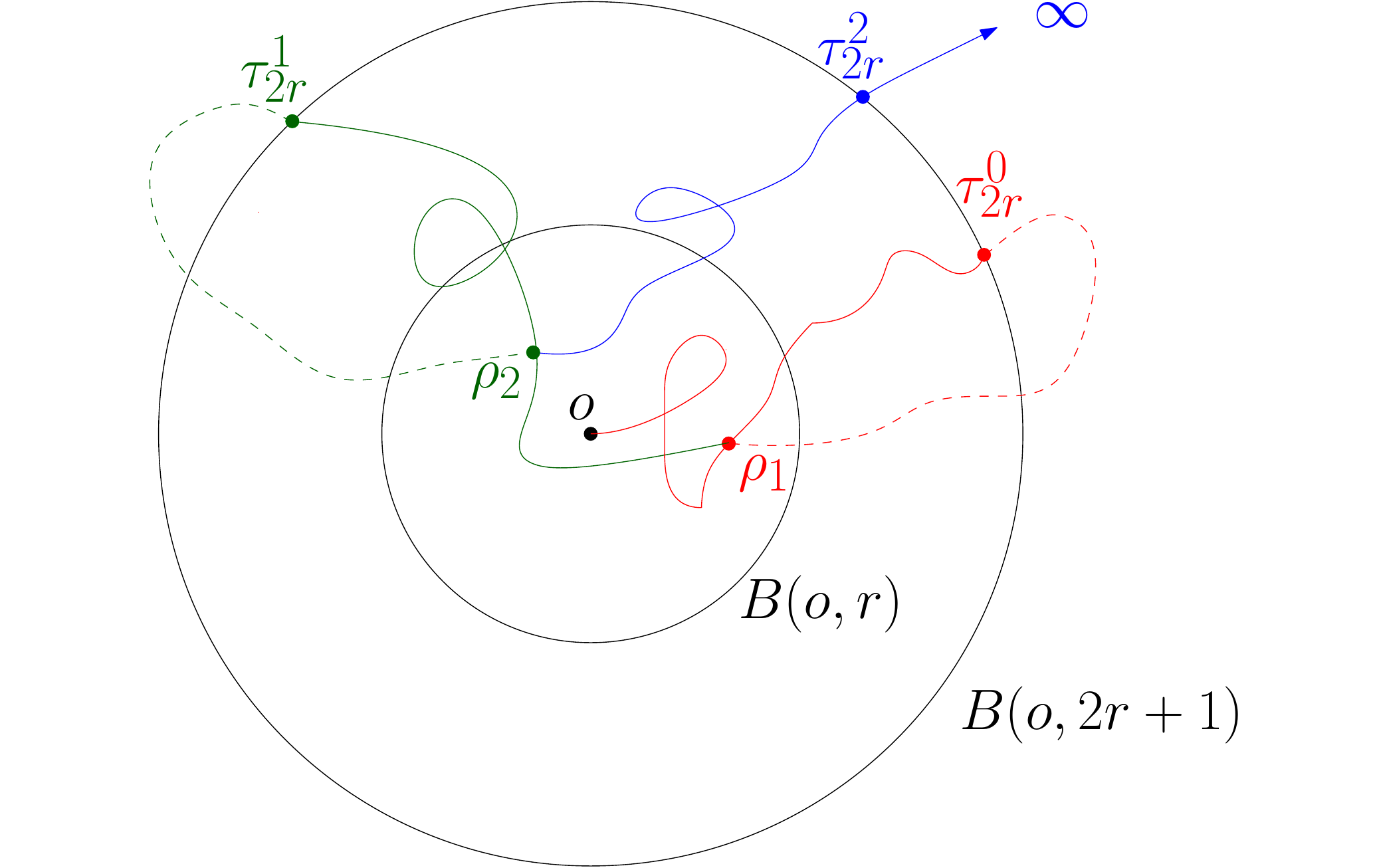} \hspace{.05\textwidth}
		\caption{Every time $S$ exits $B(o,2r)$, there is at least
			$\frac12$ chance that afterwards it never visits the vertex set in $B(o,r)$
			that is already occupied by $S$. In this figure, $\xi=3$. Note that the
			dashed part does not contribute to $\LE[S]\cap C(o,r)$. Therefore $N_r\le
			\sum_{i=0}^{\xi}(\tau^i_{2r}-\rho_i)$.  } \label{fig:geometric} \end{SCfigure}
	
	To verify this claim, we use the argument illustrated in Figure~\ref{fig:geometric}. Let $\rho_0:=0$ and $\tau^0_{2r}:=\tau_{2r}$. For $i\ge 1$, let
	\begin{align}\label{eq:time-decomposition}
	\rho_i&:=\inf\bigl\{t: t>\tau^{i-1}_{2r},\, S_t\in \LE
	[S(0,\tau^{i-1}_{2r})]\cap B(o,r)\bigr\}\\
	\noalign{and}
	\tau^{i}_{2r}&:=\inf\bigl\{t: t>\rho_i,\, S_t\notin B(o,2r)\bigr\}\,.
	\end{align}
Since $G$ has superpolynomial growth,
	by Lemma \ref{lemma:hitting-probablity}, conditioned on $S[0,\tau^{i-1}_{2r}]$,  the probability that $S$ hits a certain point in $B(o,r)$ after
	$\tau^{i-1}_{2r}$ is bounded by $cr^{-4}$, where $c$ depends only on $G$. Let $L_r$ be the occupation measure of $B(o,r)$ as defined in Theorem
	\ref{thm:occupation-time}.
	Then by conditioning on $\LE [S(0,\tau^{i-1}_{2r})]\cap B(o,r)$ and applying Theorem~\ref{thm:occupation-time}, we get
	\begin{align*}
	\P[\rho_i<\infty \mid \rho_{i-1} < \infty]
	&\le
	\Ebig{cr^{-4}\, | \LE [S(0,\tau^{i-1}_{2r})]\cap B(o,r)| \bigm| \rho_{i-1} < \infty}
	\\ &\le cr^{-4} \E[L_r] \lesssim
	r^{-1}\,.
	\end{align*}
Therefore  we may choose $r$ large enough that
	\begin{equation}\label{eq:hitting-control}
	\P[\rho_i<\infty \mid \rho_{i-1} < \infty]<1/2\,.
	\end{equation}
	Fix such an $r$.
	We have by the strong Markov property that
	\begin{equation}\label{eq:strong-Markov}
	\P[\tau^i_{2r}-\rho_i > a \mid \rho_i < \infty,\, S_{\rho_i} = x]\le
	\P[\tau_{3r} > a]
	\end{equation}
	for every $a \ge 0$ and every $x$.
	Let $\xi:=\inf\{m:\rho_m=\infty\}$. Then by \eqref{eq:hitting-control} and
	\eqref{eq:strong-Markov}, $\sum_{i=0}^{\xi-1} (\tau^i_{2r}-\rho_i)$ is
	stochastically dominated by $\sum_{i=0}^{\wt\xi-1}\tau^{i}_{3r}$, where
	$\{\tau^i_{3r} \}_{i\ge 0}$ is a sequence of i.i.d.\ random variables with
	the same distribution as $\tau_{3r}$ and $\wt\xi$ is an independent
	geometric random variable with mean $2$.
	
	Since $\LE[S]\cap C(o,r)$ is covered by the set $\bigcup_{i=0}^{\xi-1}
	S[\rho_i,\tau^i_{2r}]$ when $S(0) = o$, we have $$\E[N_r]\le \E\suml_{i=0}^{\wt
		\xi-1}\tau^{i}_{3r}=2 \E[\tau_{3r}]\,,$$ as claimed.
	
	To bound $|C(o,r)|$, we need to bound the number of vertices in $B(o,r)$
	that connect to $\trunk$ through the $\WSF$ \emph{entirely inside} $B(o,r)$.
	
	For $x,v\in B(o,r)$, write $x\simC v$ for the event that $v\in \trunk$
	and $x$ and $v$ are connected in $C(o,r)$ via a path containing no vertices
	of $\trunk$ other than $v$. For all $y\in B(o,r)$, let $\vtau_y$ be the hitting time of $y$ for a simple random walk. Let $\Psub_y$ be the distribution
	of a simple random walk $S$ starting from $y$. 
	Given $\{v_j:1\le j\le N\}\subset B(o,r)$, write $A$ for the event that $\mathrm{Ray}(o,r)=\{v_j:1\le j\le N\}$.
	For all $1\le i\le N$ and $\{v_j:1\le j\le N\}\subset B(o,r)$,
	\begin{align}
	\Pbig{y\simC v_i\bigm| A }
	&=\Psubbig_y{\textrm{$S$ hits $\trunk$ at $v_i$ and }
		\LE[S(0,\vtau_{v_i})]\subset B(o,r) \bigm| A }\nonumber\\
	\label{eq:time-reverse}
	&\le \Psubbig_{y}{\LE[S(0,\vtau_{v_i})]\subset B(o,r)}
	=\Psubbig_{v_i}{\LE[S(0,\vtau_y)]\subset B(o,r) },
	\end{align}
	where the last equality is by reversibility of LERW \cite[Lemma 7.2.1]{LawlerIntersections}.
	
	Let $M_{v}:=|\{y\in B(o,r): y\simC v\}|$. Then
	\begin{align*}
	\Ebig{M_{v_i}\bigm| A }&\le
	\suml_{y\in B(o,r)}\Psubbig_{v_i}{\LE[S(0,\vtau_y)]\subset B(o,r) }\\
	&=\Esubbig_{v_i}{\vertspace|\{y\in B(o,r): \LE[S(0,\vtau_y)]\subset B(o,r)\}|\vertspace}
	\\ &\le
	\Esubbig_{o}{\vertspace|\{y\in B(o,2r): \LE[S(0,\vtau_y)]\subset B(o,2r)\}|\vertspace}\,.
	\end{align*}
	
	Let  $\tau^i_{4r},\rho_i,\xi$ be defined  as in
	\eqref{eq:time-decomposition} but replacing $B(o,r)$ and $ B(o,2r)$ by $B(o,2r)$ and
	$B(o,4r)$, respectively. Then $\bigl\{y\in B(o,2r): \LE[S(0,\vtau_y)]\subset
	B(o,2r)\bigr\}$ is covered by the set $\bigcup_{i=0}^{\xi-1}
	S[\rho_i,\tau^i_{4r}]$ when $S(0) = o$.
	By the same argument above that proved $\E[N_r]\le 2\E[\tau_{3r}]$,  we have
	$$\Esubbig_{o}{\vertspace|\{y\in B(o,2r): \LE[S(0,\vtau_y)]\subset B(o,2r)\}|\vertspace}\le 2\E[\tau_{6r}]\,.$$
	Therefore, writing $\mathrm{Ray}(o,r)=\{v_i:1\le i\le {N_r}\}$, we have
	\begin{align*}
	\Ebig{\vertspace|C(o,r)|\vertspace}&= \Ebigg{\EBig{\suml_{i=1}^{N_r} M_{v_i}\Bigm| \mathrm{Ray}(o,r) }}\\
	&\le 2\E[\tau_{6r}]\E[N_r]\le 4\E[\tau_{6r}]^2\,. \qedhere
	\end{align*}
\end{proof}

\bibliographystyle{plain}

\end{document}